\documentclass[12pt]{amsart}
\usepackage[colorlinks=true,citecolor=black,linkcolor=black,urlcolor=blue]{hyperref}
\usepackage{amsmath,pifont,enumitem,mathtools}
\usepackage[english,  activeacute]{babel}
\usepackage[utf8]{inputenc}
\usepackage[dvipsnames]{xcolor}
\usepackage{amssymb}
\usepackage{amsthm}
\usepackage{comment}
\usepackage{graphics,graphicx}
\usepackage{array}
\usepackage{bm}
\usepackage{a4wide}
\usepackage{color, url}
\usepackage{float}
\usepackage{xcolor}

\newcommand{\seqnum}[1]{\href{http://oeis.org/#1}{\underline{#1}}}

\theoremstyle{plain}
\newtheorem{theorem}{Theorem}[section]

\newtheorem{coro}[theorem]{Corollary}

\theoremstyle{definition}

\DeclareRobustCommand{\stirling}{\genfrac\{\}{0pt}{}}
\newcommand{\Sn}{{\mathfrak S}}

\newcommand{\I}{{\mathcal I}}
\newcommand{\B}{{\mathcal B}}
\newcommand{\F}{{\mathcal F}}
\newcommand{\V}{{\mathcal V}}
\newcommand{\id}{{\mathcal{ID}}}
\def\val{{\textsf{val}}}
\def\run{{\textsf{run}}}
\def\inc{{\textsf{inc}}}
\def\tot{{\textsf{Tot}}}

\def\des{{\textsf{des}}}
\def\asc{{\textsf{asc}}}
\def\rlm{{\textsf{rlm}}}
\def\peak{{\textsf{peak}}}
\def\s{{\textsf{s}}}
\def\rp{{\textsf{rp}}}
\def\dasc{{\textsf{dasc}}}

\newcommand{\jluc}[1]{\mbox{}{\sf\color{blue}[J-Luc: #1]}\marginpar{\color{magenta}\Large$*$}}

\title{Some distributions in increasing and flattened permutations}

\date{\today}

\begin{document}

\author[J.-L Baril]{Jean-Luc Baril}
\address{LIB, Universit\'e de Bourgogne Franche-Comt\'e,   B.P. 47 870, 21078, Dijon Cedex, France}
\email{barjl@u-bourgogne.fr}

\author[J. L. Ram\'{\i}rez]{Jos\'e L. Ram\'{\i}rez}
\address{Departamento de Matem\'aticas,  Universidad Nacional de Colombia,  Bogot\'a, Colombia}
\email{jlramirezr@unal.edu.co}

\begin{abstract} We examine  the distribution and popularity of different parameters (such as the number of descents, runs, valleys, peaks, right-to-left minima, and more) on the sets of increasing and flattened permutations. For each parameter, we provide an exponential generating function for its corresponding distribution and popularity.Additionally, we present one-to-one correspondences between these permutations and some classes of simpler combinatorial objects.
\end{abstract}

\maketitle

\section{Introduction and Notation}

Let $\Sn_n$ denote  the set of permutations of $[n]:=\{1,2,\ldots,n\}$. For convenience, we define $\Sn_0$ to consist of the empty permutation $\epsilon$, and  we set $\Sn=\bigcup_{n\geq 0}\Sn_n$. A permutation $\pi\in\Sn$ can be written in \emph{one-line notation} as the word $\pi = \pi_1\pi_2\cdots \pi_n$, where $\pi_i=\pi(i)$ for $1\leq i \leq n$, and it is represented by the empty word (also denoted $\epsilon$) if $\pi\in \Sn_0$. We denote the length of $\pi$ by $|\pi|$. Given a sequence of distinct integers $a=a_1a_2\cdots a_n$, $\texttt{red}(a)$ is the unique permutation in $\Sn_n$ order-isomorphic to $a$, and $\texttt{dom}(a)$ is the set $\{a_1,a_2,\ldots , a_n\}$. For instance, if $a=3~6~1~4$ then $\texttt{red}(a)=2~4~1~3$ and $\texttt{dom}(a)=\{1,3,4,6\}$.

Let $\pi=\pi_1\pi_2\cdots\pi_n\in \Sn_n$.  A \emph{valley} in $\pi$ is an index $\ell$, where  $2\leq \ell \leq n-1$, such that $\pi_{\ell-1}>\pi_\ell<\pi_{\ell+1}$.  The \emph{height} of a valley is given by the image of its index, that is $\pi_\ell$. The number of valleys in $\pi$ is denoted by $\val(\pi)$. Analogously, we
define the peak statistic.   A \emph{peak} in $\pi$ is an index $\ell$, where  $2\leq \ell \leq n-1$, such that $\pi_{\ell-1}<\pi_\ell>\pi_{\ell+1}$.  The number of peaks in $\pi$ is denoted by $\peak(\pi)$.   We say that $\pi$ has an \emph{ascent} (resp. \emph{descent}) at position $\ell$ if $\pi_\ell < \pi_{\ell+1}$ (resp. $\pi_\ell > \pi_{\ell+1}$), where $\ell \in [n]$.  The number of ascents and descents of $\pi$ is denoted by $\asc(\pi)$ and $\des(\pi)$, respectively. Obviously, if $\pi\in\Sn_n$ then we have $\des(\pi)=n-1-\asc(\pi)$. A \emph{run} in a permutation $\pi$ is a subword $\pi_\ell\pi_{\ell+1}\cdots \pi_{\ell+s} \pi_{\ell+s+1}$, where $\ell, \ell+1,\dots, \ell+s$ are consecutive ascents. The number of runs of $\pi$ is denoted by $\run(\pi)$.  A \emph{right-to-left minimum} of $\pi$ is an entry $\pi_i$ such that if $j>i$, then  $\pi_i< \pi_j$. The number of right-to-left minima of $\pi$ is denoted by $\rlm(\pi)$. All these parameters are summarized in Table \ref{parameters}.

\begin{table}[htp!]
\begin{tabular}{|ll|ll|}\hline
$\val(\pi)$  &  number of valleys of $\pi$ 		& $\run(\pi)$   &  number of runs of $\pi$\\ \hline
$\des(\pi)$ & number of descents of $\pi$   &  $\asc(\pi)$ & number of ascents of $\pi$ \\  \hline $\rlm(\pi)$  & number of right-to-left minima of $\pi$	&		
$\peak(\pi)$ & number of peaks of $\pi$ 	\\ \hline
\end{tabular}\vspace{4mm}
\caption{Parameters.}\label{parameters}
\end{table}

In this paper, we focus on the distribution of the aforementioned permutation statistics within two families of permutations: {\it increasing} and {\it flattened} permutations. 

\begin{figure}[ht!]
\centering
\includegraphics[scale=0.9]{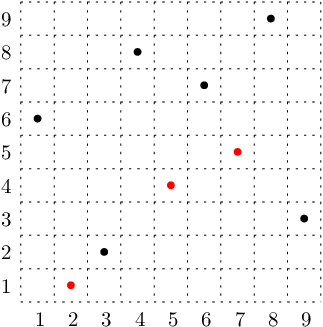}\qquad\qquad\includegraphics[scale=0.9]{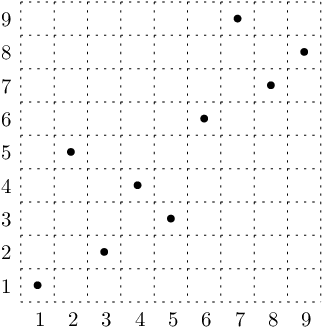}
\caption{Increasing permutation $\pi=6~1~2~8~4~7~5~9~3$, and flattened permutation $\pi=1~5~2~4~3~6~9~7~8$.} \label{fig1}
\end{figure}

A permutation $\pi$ is called \emph{increasing} if the heights of its valleys form an increasing sequence from left to right. For example,  in Figure~\ref{fig1}, we illustrate an increasing permutation of length $9$, where red vertices mark  its valleys.   We denote by $\I$ the set of all increasing permutations, and by $\I_n$ the set of those of length $n$. This concept has also been studied in the context of other combinatorial structures. One of the earliest works was by Barcucci et al.\ \cite{barcucci} in the context of Dyck path. For example,  they prove that the number of Dyck path of length $2n$, where the heights of their valleys form a non-decreasing sequence, is given by the Fibonacci number $F_{2n-1}$ (see also \cite{CFJR, FR}). This concept was later extended by Fl\'orez and Ram\'irez \cite{FR2} for Motzkin paths. Fl\'orez at al.\ \cite{FRV} also studied this concept in the context of integer compositions, showing connections with unimodal compositions.

A permutation $\pi$ is called \emph{flattened}  if it consists of runs arranged from left to right, with the first entries of each run in increasing order.  We denote by $\F$ the set of all flattened permutations and by $\F_n$ the set of those of length $n$, and we clearly have $\F_n\subset\I_n$ for $n\geq 1$.  For example,  the permutation $\pi=1~5~2~4~3~6~9~7~8$ is flattened, with runs $1~5, 2~4, 3~6~9, 7~8$. This concept was introduced by Callan \cite{Callan}, who determined the number of flattened partitions of length $n$ that avoid a pattern of length  3. Since then, several authors have investigated this concept. For instance,  Nabawanda et al.\ \cite{naraba} studied the run distribution of flattened permutations and provided a bijection between flattened partitions over $[n+1]$ and set partitions of $[n]$. Mansour et al.\ \cite{manshawag,manshawang} examined the avoidance  of subword patterns in flattened partitions  (see also \cite{ manss}) and the distribution of consecutive patterns in flattened permutations (see \cite{msha}).  Notably, Mansour et al.\ \cite{manshawag} also studied some statistics such as descents, ascents, peaks, and valleys in a different context of flattened permutations, specifically focusing on the cycle notation of permutations. Recently, this concept has been explored in other combinatorial objects. For example, Buck et al.\ \cite{FlatStirling}  introduced the notion of Stirling  permutations, Elder et al.\ \cite{flat_pf} defined and enumerated flattened parking functions, and  Baril et al.\ \cite{BHR} investigated several statistics on flattened  Catalan words.

Permutations are labeled combinatorial objects, and so
it is appropriate to use exponential generating functions (e.g.f.) for the enumerative analysis of them. We provide bivariate exponential generating function for increasing and flattened permutations with respect to the length and the parameters listed in Table~\ref{parameters}, that is, function where the coefficient of $\frac{x^ny^k}{n!}$ in its series expansion is the number of permutations (increasing or flattened) of length $n$ such that the parameter equals $k$. We also give the e.g.f. for the popularity of the parameter,
function where the coefficient of  $\frac{x^n}{n!}$ in its series expansion is the total number of this parameter appearing in all
length $n$ permutations class. 

Our method consists in constructing recursively the combinatorial class in question
from  smaller classes, $\mathcal{A}_1$ and $\mathcal{A}_2$,
using the usual labeled product $\mathcal{A}_1\star\mathcal{A}_2$ and the boxed product
$\mathcal{A}_1^\square\star\mathcal{A}_2$.
The boxed product $\mathcal{A}_1^\square\star\mathcal{A}_2$ is a subset of $\mathcal{A}_1\star\mathcal{A}_2$ where the smallest label appears in the $\mathcal{A}_1$ component. For more details about labeled combinatorial structures, the boxed product, and their applications, refer to the book of Flajolet and  Sedgewick \cite{Flajolet2}, as well as \cite{bar} for several applications.

\section{Enumeration of Increasing Permutations}
In this section, we are interested in the bivariate  generating function
\begin{align}
I^s(x,y):=\sum_{\pi \in\I}\frac{x^{|\pi|}}{|\pi|!}y^{\s(\pi)}=\sum_{n\geq 0}\frac{x^{n}}{n!}\sum_{\pi\in\I_n}y^{\s(\pi)}.
\end{align} whenever the statistic $\s$ is $\val$, $\run=1+\des$, and $\rlm$.
The coefficient of $\frac{x^ny^k}{n!}$ in $I^s(x,y)$ is the number of increasing permutations $\pi$ of  length $n$ satisfying $s(\pi)=k$.

\subsection{For the number of valleys: $\s=\val$} In this part, we focus on the enumeration of increasing permutations with respect to the length and the statistic $\s=\val$, which counts the number of valleys.

For the permutations, Elizalde and Noy \cite{ENoy} established that  the exponential generating function $P(x,y)$ for the number of permutations with respect to the number of valleys and the length is given by
\begin{align}
P(x,y)=\sum_{n\geq 0}\frac{x^{n}}{n!}\sum_{\pi\in\Sn_n}y^{\val(\pi)}=\frac{\sqrt{1-y}}{\sqrt{1-y} - \tanh(x\sqrt{1-y})} \label{gfVall}.
\end{align}
The coefficient of $\frac{x^ny^k}{n!}$ in the series expansion of $P(x,y)$ corresponds to the array \seqnum{A008303}  in the OEIS \cite{OEIS}.  This counting problem has been extensively studied, with many related papers addressing various aspects of it, see for example \cite{Bukata, MR, Z}.

We denote by $\V$ the set of \emph{valleyless permutations}, which are permutations free of valleys. Rieper and Zeleke \cite{RZ} established a bijection between valleyless permutations and integer compositions.  Consequently, the exponential generating function for valleyless permutations is given by
 $$V(x)=\frac{1}{2} \left(e^{2 x}+1\right)=1+e^x\sinh x.$$

\begin{theorem}\label{teo1}
The generating function for increasing permutations with respect to the
length and the number of valleys is
$$I^\val(x,y)=e^{2 x - \frac{y}{4} + \frac{1}{4} e^{2 x}y - \frac{x y}{2}}\left(1+\frac{e^{y/4}}{2}\int_0^xe^{\frac{-e^{2t}}{4}y - 2t + \frac{yt}{2}}(-2 + y - e^{2 t} y)dt\right).$$
\end{theorem}
\begin{proof}
Let $\pi$ be a nonempty increasing permutation. We consider the decomposition of $\pi$ based on the position of the entry 1.  Specifically, $\pi$ can be expressed in one of the following forms:  $\pi=\pi_11$,  $\pi=1\pi_2$, or $\pi=\pi'1\pi_2$, where $\pi_2\neq \epsilon$ and $\pi'$ is free of valleys. From this decomposition, we derive the following symbolic equation:
$$\I=\epsilon + \I \star \{x\}^{\square} + \{x\}^{\square}\star(\I-\epsilon)  + (\V-\epsilon) \star \{x\}^{\square} \star(\I-\epsilon).$$
Here, $\mathcal{A}^{\square} \star \mathcal{B}$ denotes the boxed product of the combinatorial classes $\mathcal{A}$ and $\mathcal{B}$.   From the symbolic equation,  and after multiplying by $y$ whenever a valley is created, we obtain the functional equation
$$I^\val(x,y)= 1 +  \int_{0}^x I^\val(t,y)dt + \int_{0}^x (I^\val(t,y)-1)dt + y\int_{0}^x e^t \sinh t (I^\val(t,y)-1)dt.$$
Differentiating this functional  equation yields the following differential equation:
$$\frac{d}{dx} I^\val(x,y)=    2I^\val(x,y) -1 + y e^x \sinh x (I^\val(x,y)-1), \quad I^\val(0,y)=1.$$
The solution of this differential equation provides the desired result. 
\end{proof}

Notice that  the solution $I^\val(x,y)$ can be expressed  in a closed form without integrals, depending on the WhittakerM function \cite{whi}. However, this form is not particularly elegant, so we do not include it here. We solve the differential equations of this paper by using symbolic software. 

The first few terms of the formal power series are
\begin{multline*}
I^\val(x,y)=\frac{1}{2} \left(e^{2 x}+1\right)+\frac{1}{8} \left(-4 e^{2 x} x+e^{4 x}-1\right)
   y\\
   +\frac{1}{64} \left(e^{2 x} \left(8 x^2-4 x-1\right)+e^{4 x} (2-8 x)+e^{6 x}-2\right)
   y^2+ O(y^3).
   \end{multline*}

Let $\inc^\val(n,k)$ denote the number of increasing permutations in $\Sn_n$ with exactly $k$ valleys. The first few values of this sequence are \footnotesize
$$I^\val:=[\inc^\val(n,k)]_{n, k \geq 0}= \left(
\begin{array}{cccccc}
 1 & 0 & 0 & 0 & 0 &\cdots\\
 1 & 0 & 0 & 0 & 0&\cdots \\
 2 & 0 & 0 & 0 & 0 &\cdots\\
 4 & 2 & 0 & 0 & 0 &\cdots\\
 8 &  \framebox{\textbf{16}} & 0 & 0 & 0 &\cdots\\
 16 & 88 & 8 & 0 & 0&\cdots\\
 32 & 416 & 136 & 0 & 0 &\cdots\\
 64 & 1824 & 1440 & 48 & 0 &\cdots\\
 128 & 7680 & 12288 & 1384 & 0 &\cdots\\
 \vdots &\vdots & \vdots & \vdots & \vdots &\ddots\\
\end{array}
\right).$$
\normalsize 
For example, $\inc^\val(4,1)=16$, which corresponds to the  boxed entry in the array $I^{\val}$ above. The 16 increasing permutations corresponding to this value are
\begin{align*}
&1324, \quad 1423, \quad 2134, \quad 2143, \quad 2314, \quad 2413, \quad 3124, \quad 3142,\\
&3214,  \quad 3241, \quad 3412,  \quad 4123, \quad  4132,  \quad 4213,  \quad  4231, \quad 4312.
\end{align*}
Note that $\inc^\val(n,1)=2^{n-2}(2^{n-1} - n)$, for $n\geq 3$, which corresponds to the sequence \seqnum{A000431} in \cite{OEIS}.

Due to the functional equation involving $I^\val(x,y)$, and after observing that  
$$\V(x)=1+\sum\limits_{n\geq 1}2^{n-1}\frac{x^n}{n!},$$ a straightforward series calculation leads to the following result. 

\begin{coro} For $n,k\geq 1$, we have the recursive formula
$$\inc^\val(n+1,k)=2\, \inc^\val(n,k)+\sum\limits_{i=1}^{n-1}{n \choose i}2^{i-1}\inc^\val(n-i,k-1),$$
with the initial condition $\inc^\val(0,0)=1$ and  $\inc^\val(n,0)=2^{n-1}$ for $n\geq 1$.
\end{coro}

\begin{coro}\label{corotoinper}
The generating function $I(x)$ for total number of increasing  permutations is
\begin{align*}
I(x):=I^\val(x,1)&=e^{\frac{3 x}{2} - \frac{1}{4} + \frac{1}{4} e^{2 x}}\left(1-\frac{e^{1/4}}{2}\int_0^xe^{\frac{-e^{2t}}{4} -  \frac{3t}{2}}(1 + e^{2 t})dt\right)\\
&=1  + \frac{x}{1!}  + 2\frac{x^2}{2!}  + 6\frac{x^3}{3!} +24 \frac{x^4}{4!}  + 112\frac{x^5}{5!}  + 584\frac{x^6}{6!} ++ 3376\frac{x^7}{7!}+ O(x^8).
\end{align*}
Moreover, if $\inc(n)$ is the number of increasing permutations of length $n$, then 
$$\inc(n+1)=2\,\inc(n)+\sum\limits_{i=1}^{n-1}{n \choose i}2^{i-1}\inc(n-i), \quad n\geq 1$$
with the initial conditions $\inc(0)=1=\inc(1)$.
\end{coro}

\begin{coro}
The generating function for total number of valleys in all increasing  permutations of a given length is
\begin{align*}
\tot^\val(x)&=\partial_y(I^\val(x,y))\rvert_{y=1} \\
&= 2\frac{x^3}{3!} +16 \frac{x^4}{4!}  + 104\frac{x^5}{5!}  + 688\frac{x^6}{6!} + 4848\frac{x^7}{7!}+ O(x^8).
\end{align*}
\end{coro}
The sequences of coefficients of $\frac{x^n}{n!}$ from the two series expansions in the previous corollaries do not appear in the OEIS.

\subsection{For the number of runs: \s=\run=1+\des=n-\asc }
In this part, we focus on the enumeration of increasing permutations with respect to the length and the statistic $\s=\run$ that counts the number of runs. It is worth noticing that, in a nonempty permutation, the number of descents equals the number of runs minus one, i.e.,  $\des=\run-1$. Therefore, we have $$I^\run(x,y)=1+y(I^\des(x,y)-1)\quad\mbox{ and } \quad I^\asc(x,y)=1+\frac{1}{y}(I^\des(xy,1/y)-1).$$

For the classical  permutations, the number of permutations of length $n$ with exactly $k$ descents is given by the Eulerian numbers (sequence \seqnum{A008292}). It is well-known that (cf. \cite{bar,BookMezo,peter})  
\begin{align*}
E(x,y)=\sum_{\pi\in \Sn}\frac{x^{|\pi|}}{|\pi|!}y^{\des(\pi)}=\frac{1-y}{1-ye^{x(1-y)}}.
\end{align*}

\begin{theorem}
The generating function for increasing permutations with respect to the
number of runs is

\resizebox{6.4in}{!}{$I^\run(x,y)=\exp \left(\frac{x (y+1) \left(y^2+y+1\right)+y \left(e^{x
   (y+1)}-1\right)}{(y+1)^2}\right) \left(\cfrac{e^{y/(y+1)^2}}{y+1}\displaystyle \int_0^x
   \left( -ye^{(y+1)t}-1\right) \exp \left(-\frac{(y+1)
   \left(y^2+y+1\right)t+y e^{(y+1)t}}{(y+1)^2}\right) \,
   dt+1\right).$ }

\end{theorem}
\begin{proof}
Let $\pi$ be a nonempty increasing permutation. We consider the same decomposition of $\pi$ based on the position of $1$, as in the proof of Theorem~\ref{teo1}. From this decomposition, we obtain the symbolic equation:
\begin{align}\label{ec1}
\I=\epsilon + \{x\}+(\I-\epsilon) \star \{x\}^{\square} + \{x\}^{\square}\star (\I-\epsilon)  + (\V-\epsilon) \star \{x\}^{\square} \star(\I-\epsilon),
\end{align}
where $\V$ is the set of valleyless permutations.

Now let us give the exponential generating function for the set $\V$ with respect to the length and the number of runs. Using the following recursive decomposition of $\V$, 
$$\V=\epsilon + \{x\}+\{x\}^{\square}\star (\V-\epsilon)+(\V-\epsilon) \star \{x\}^{\square},$$ 
we deduce the functional equation (we multiply by $y$ whenever a run is created)
$$V(x,y)=1+xy+\int_{0}^x (V(t,y)-1)dt+y\int_{0}^x (V(t,y)-1)dt.$$
Then we obtain $$V(x,y)=\frac{e^{\left(y +1\right) x} y +1}{y +1}=1+\sum_{n\geq 1}\sum_{j=0}^{n-1}\binom{n-1}{j}y^{j+1}\frac{x^n}{n!}.$$
From \eqref{ec1}, and after multiplying by $y$ whenever a run is created, we obtain the functional equation
$$I^\run(x,y)= 1 +xy +(1+y) \int_{0}^x (I^\run(t,y)-1)dt  + \int_{0}^x (V(t,y)-1)(I^\run(t,y)-1)dt.$$
Differentiating the previous equation, we obtain the following differential equation:
$$\frac{d}{dx} I^\run(x,y)=  y+ (1+y)(I^\run(x,y) -1) + (V(x,y)-1)(I^\run(x,y)-1), \quad I^\run(0,y)=1.$$
The solution of this differential equation is the desired result. 
\end{proof}

The first few terms of
the formal power series are
\begin{multline*}
I^\run(x,y)=1+\left(e^x-1\right) y+e^x \left(-x+e^x-1\right) y^2 + \\\frac{1}{2} \left(e^x \left(-x^2+4 x+1\right)-4
   e^{2 x}+e^{3 x}+2\right) y^3+O\left(y^4\right).
   \end{multline*}
Let $\inc^\run(n,k)$ the number of increasing permutations in $\Sn_n$ with exactly $k$ runs. 
Then, we have the following \footnotesize
$$I^\run:=[\inc^\run(n,k)]_{n, k \geq 0}=\left(
\begin{array}{cccccccc}
 1 & 0 & 0 & 0 & 0 &0&0&\cdots\\
 0 & 1 & 0 & 0 & 0&0&0&\cdots \\
 0 & 1 & 1 & 0 & 0 &0&0&\cdots\\
 0 & 1 & 4 & 1 & 0 &0&0&\cdots\\
 0 & 1 & \framebox{\textbf{11}} & 11 & 1 &0&0&\cdots\\
 0 & 1& 26 & 58 & 26 & 1&0&\cdots\\
 0 & 1& 57 & 234 & 234 &57&1& \cdots\\
  0 & 1& 120 & 831 & 1472 &831&120& \cdots\\
   0 & 1& 247 & 2757 & 7735 &7735&2757& \cdots\\
 \vdots &\vdots & \vdots & \vdots & \vdots& \vdots & \vdots&\ddots\\
\end{array}
\right).$$ \normalsize
For example, $\inc^\run(4,2)=11$, the entry boxed in the array $I^\run$ above, and the corresponding increasing permutations are
\begin{align}\label{ex1}
&1243, \quad 1324, \quad 1342, \quad 1423, \quad 2134, \quad 2314, \\
&2341, \quad 2413, \quad 3124,  \quad 3412, \quad 4123. \nonumber
\end{align}

Note that the sequence $\inc^\run(n,2)$, $n\geq 2$,  coincides with the sequence of  Eulerian numbers \seqnum{A000295} that enumerate permutations with only one descent (which are exactly  increasing permutations with two runs). 

A simple calculation on series induces the following result.  
\begin{coro} We have the following recursive formula for $n,k\geq 1$
$$\inc^\run(n+1,k)= \inc^\run(n,k)+\inc^\run(n,k-1)+\sum\limits_{i=1}^{n-1}\sum\limits_{j=1}^{k-1}{n \choose i}{i-1\choose k-j-1}\inc^\run(n-i,j),$$
with the initial condition $\inc^\run(0,0)=1$, $\inc^\run(n,0)=0$ for $n\geq 1$.
\end{coro}

Observing the array $I^\run$, we can hope the existence of an involution $\phi$ on $\I_n$  such that $\run(\phi(\pi))=n+1-\run(\pi)$. Then, let us define recursively the map $\phi$ on $\I$ as follows:

For $\pi\in\I$, we set:
$$
\phi(\pi)=\left\{\begin{array}{llr}
\epsilon&\text{if }\pi=\epsilon&(i)\\
(1+\phi(\pi'))1&\text{if }\pi=1(1+\pi'), \quad \pi'\in \I &(ii)\\
1(1+\phi(\pi'))&\text{if }\pi=(1+\pi')1, \quad\pi'\in \I, \pi'\neq\epsilon&(iii)\\
\pi_1^r1\chi(\phi(\texttt{red}(\pi_2)),\texttt{dom}(\pi_2))&\text{if }\pi=\pi_11\pi_2, \quad \pi_1,\pi_2\neq \epsilon&(iv)\\\end{array}\right. $$
where $\pi_1^r$ is the reverse of $\pi_1$ and $\texttt{dom}(\pi_2)$ is the set of values in $\pi_2$ and $\chi(\pi,S)$ (with $|\pi|=|S|$) is the unique sequence  taking values in $S$ such that $\texttt{red}(\sigma)=\pi$.
For instance, if $\pi=4132$ and $S=\{2,6,7,9\}$,  then $\chi(\pi,S)=9276$ since $\texttt{red}(9276)=4132$.

Let us give an example. Easily, we can see that $\phi(1)=1$, $\phi(12)=21$ and $\phi(21)=12$. Now, let $\pi=23154$ be an increasing permutation with 3 runs. Then, we have  $\pi=\pi_11\pi2$ with $\pi_1=23$, $\pi_2=54$ and $\texttt{dom}(\pi_2)=\{4,5\}$. Therefore, we obtain 
$$\chi(\phi(\texttt{red}(\pi_2),\texttt{dom}(\pi_2))=\chi(\phi(21),\{4,5\})=\chi(12,\{4,5\})=45.$$
Finally, we have $\phi(23154)=32145$ and the number of runs in $32145$ is $3$ which also is $5+1$ minus the number of runs of $\pi=23154$.

For example, the  images by $\phi$ (keeping the order) of the permutations  in \eqref{ex1} are the increasing permutations with 3 runs, respectively:
\begin{align*}
&3421, \quad 3241, \quad 2431, \quad 4231, \quad 2143, \quad 3214,\\
&1432, \quad 4213, \quad 3142,  \quad 4312, \quad 4132.   
\end{align*}

\begin{coro}
The generating function for total number of runs in all increasing  permutations of a given length is
\begin{align*}
\tot^\run(x)&=\partial_y(I^\run(x,y))\rvert_{y=1} \\
&= x+3\frac{x^2}{2!}+12\frac{x^3}{3!} +60 \frac{x^4}{4!}  + 336\frac{x^5}{5!}  + 2044\frac{x^6}{6!} + 13504\frac{x^7}{7!}+ O(x^8).
\end{align*}
\end{coro}
The sequence of coefficients of $\frac{x^n}{n!}$ in this series expansions does not appear in the OEIS.
\medskip

Let $\I^{\val=\des}_n$ be the set of increasing permutations $\pi$ of length $n$ such that   $\val(\pi)=\des(\pi)$. We set $\I^{\val=\des}=\bigcup_{n\geq 0}\I^{\val=\des}_n$.

\begin{theorem}  The cardinality $\inc^{\val=\des}(n)$ of $\I^{\val=\des}_n$ is given by the $n$-th Gould number (see \seqnum{A040027} in \cite{OEIS}):
$$\inc^{\val=\des}(n)=\sum\limits_{k=0}^{n-1}\sum\limits_{j=k}^{n-1}\stirling{j}{k}k^{n-1-j},$$
where $\stirling{j}{k}$ is the Stirling number of second kind.
\end{theorem}
\begin{proof} Let $\pi$ be a permutation in $\I^{\val=\des}_n$. Then $\pi$ has a unique decomposition as follows: $\pi=\pi_1a_1\pi_2a_2 \cdots \pi_ka_k$, where $a_1, a_2, \ldots, a_k$ are the right-to-left minima of $\pi$. Due to the equality $\val(\pi)=\des(\pi)$, $\pi_ia_i$ is either reduced to $a_i$ whenever $\pi_i=\epsilon$, or $\pi_ia_i$ has only one descent and this descent involves $a_i$. Moreover, $\pi_k$ is necessarily empty. Therefore, we construct the partition $B=B_1/B_2/\cdots /B_k$ of $[n]$  where $B_i=\texttt{dom}(\pi_ia_i)$ for $1\leq i\leq k$. Since the $a_i$ are right-to-left minima of $\pi$, then we have $a_1=1<a_2<\cdots<a_k$, and then $\min(B_1)<\min(B_2)<\cdots<\min(B_k)=a_k$. Thus, $B$ is a partition ordered by the smallest element of the blocks such that the last block is reduced to a singleton.

Conversely, let us consider such a partition $B=B_1/B_2/\cdots/B_k$. Then we construct an permutation $\pi$ as follows: $\pi=\pi_1a_1\pi_2a_2 \cdots \pi_ka_k$, where $a_i=\min(B_i)$ and $\pi_i$ consists of elements of $B_i\backslash\{a_i\}$ sorted in increasing order. A simple observation proves that $\pi\in \I^{\val=\des}_n$, which establishes a bijection between $\I^{\val=\des}_n$ and the set of partitions of $[n]$ having the last part reduced to a singleton. Since these partitions are counted by the $n$th Gould number (see \seqnum{A040027} in \cite{OEIS}), we obtain the desired result.
\end{proof}

\subsection{For the number of right-to-left minima: \s=\rlm}

In this part, we focus on the enumeration of increasing permutations with respect to the length and the statistic $\s=\rlm$ that counts the number of right-to-left minima. 

For the classical  permutations, the number of permutations of length $n$ with exactly $k$ right-to-left minima is the (unsigned) Stirling  numbers of the first kind (sequence \seqnum{A132393}). It is well-known that (cf. \cite{BookMezo}) the corresponding bivariate generating function is  
\begin{align*}
S(x,y)=\sum_{\pi\in \Sn}\frac{x^{|\pi|}}{|\pi|!}y^{\rlm(\pi)}=\frac{1}{(1-x)^y}.
\end{align*}

\begin{theorem}
The generating function for increasing permutations with respect to the
number of right-to-left minima is
$$I^\rlm(x,y)={e^{\frac{y}{4}\left( 2\,x+e^{2x} \right)}} \left( \frac{y}{2}\int _{0}^{x}\! \left( 2\,I \left( t \right) -1-{e^{2t}
} \right) {e^{-\frac{y}{4} \left( 2 t+{e^{2t}} \right) }}{dt
}+{e^{-\frac{y}{4}}} \right),
$$ where $I(z)$ is given in Corollary \ref{corotoinper}.
\end{theorem}
\begin{proof}
Let $\pi$ be a nonempty increasing permutation. We consider the decomposition of $\pi$ respect to the position of the entry 1.  From this decomposition we obtain the symbolic equation:
$$\I=\epsilon + \I \star \{x\}^{\square} + \{x\}^{\square}\star(\I-\epsilon)  + (\V-\epsilon) \star \{x\}^{\square} \star(\I-\epsilon).$$
From the symbolic method (cf. \cite{Flajolet2}) we obtain the functional equation
$$I^\rlm(x,y)= 1 +  y\int_{0}^x I^\rlm(t,1)dt + y\int_{0}^x (I^\rlm(t,y)-1)dt + y\int_{0}^x e^t \sinh t (I^\rlm(t,y)-1)dt,$$
where $I^\rlm(t,1)=I(t)$ is the exponential generating function for increasing permutations with respect to the length (see Corollary~\ref{corotoinper}).

Differentiating the previous equation, we obtain the following differential equation:
$$\frac{d}{dx} I^\rlm(x,y)=    yI^\rlm(x,1) + yI^\rlm(x,y)-y+ y e^x \sinh x (I^\rlm(x,y)-1), \quad I^\rlm(0,y)=1.$$
The solution of this equation is the desired result. 
\end{proof}

Let $\inc^\rlm(n,k)$ the number of increasing permutations in $\Sn_n$ with exactly $k$ right-to-left minima. 
Then, we have the following  \footnotesize
$$I^\rlm:=[\inc^\rlm(n,k)]_{n, k \geq 0}=\left(
\begin{array}{cccccccc}
 1 & 0 & 0 & 0 & 0 &0&0&\cdots\\
 0 & 1 & 0 & 0 & 0&0&0&\cdots \\
 0 & 1 & 1 & 0 & 0 &0&0&\cdots\\
 0 & 2 & 3 & 1 & 0 &0&0&\cdots\\
 0 & 6 & \framebox{\textbf{11}} & 6 & 1 &0&0&\cdots\\
 0 & 24& 42 & 35 & 10 & 1&0&\cdots\\
 0 & 112& 174 & 197 & 85 &15&1& \cdots\\
  0 & 584& 812 & 1116 & 667 &175&21& \cdots\\
   0 & 3376& 4836 & 6510 & 5085 &1822&322& \cdots\\
 \vdots &\vdots & \vdots & \vdots & \vdots& \vdots & \vdots&\ddots\\
\end{array}
\right).$$ \normalsize

For example, $\inc^\rlm(4,2)=11$, the entry boxed in the array $I^\rlm$ above, and the corresponding increasing permutations are
\begin{align*}
&1342, \quad 1432, \quad 2143, \quad 2314, \quad 2413, \quad 3142,\\
&3214, \quad 3412, \quad 4132,  \quad 4213, \quad 4312.
\end{align*}

Due to the functional equation involving $I^\rlm(x,y)$ a simple calculation on series induces the following corollary. 

\begin{coro} We have the following recursive formula for $n,k\geq 1$:
$$\inc^\rlm(n+1,k)=i^\rlm(n,k-1)+\sum\limits_{i=1}^{n-1}{n \choose i}2^{i-1}\inc^\rlm(n-i,k-1),$$
with the initial condition $\inc^\rlm(0,0)=1$, $\inc^\rlm(n,0)=0$ for $n\geq 1$, and  $\inc^\rlm(n,1)=\inc(n-1)$ for $n\geq 1$.
\end{coro}
\begin{coro}
The generating function for total number of right-to-left minima in all increasing  permutations of a given length is
\begin{align*}
\tot^\rlm(x)&=\partial_y(I^\rlm(x,y))\rvert_{y=1} \\
&=  x+ 3\frac{x^2}{2!}+11\frac{x^3}{3!} +50\frac{x^4}{4!}  + 258\frac{x^5}{5!}  + 1472\frac{x^6}{6!}+9232\frac{x^7}{7!}+O(x^8).
\end{align*}
\end{coro}
The sequence of coefficients of $\frac{x^n}{n!}$ in this series expansions does not appear in \cite{OEIS}.

\subsection{For the number of peaks: \s=peak}
In this part, we focus on the enumeration of increasing permutations with respect to the length and the statistic $\s=\peak$ that counts the number of peaks. Let $I^\peak(x,y)$ be the bivariate generating function for this. We will consider four sets $\mathcal{I}^{a,b}$, $a,b\in\{0,1\}$, consisting of increasing permutations $\pi$ such that $\pi$ does not start with a descent if and only if $a=0$, and does not end with an ascent if and only if $b=0$. For $a,b\in\{0,1\}$,  $I^{a,b}(x,y)$ is the generating function for the number of permutations $\pi\in\mathcal{I}^{a,b}$ with respect to the length and the number of valleys.

\begin{theorem} The generating function for increasing permutations with respect to the number of peaks is $$I^\peak(x,y)=1+x+y(I^{0,0}(x,y)-1-x)+I^{0,1}(x,y)+\frac{1}{y}I^{1,1}(x,y)+I^{1,0}(x,y),$$
where    
\begin{align*}
I^{0,1}(x,y)&=e^{\frac{y e^{2 x}}{4}-\frac{x y}{2}+x}\int_{0}^{x}e^{-\frac{y e^{2 t}}{4}+\frac{t y}{2}-t}dt \,  -\\ &\qquad\qquad\qquad y \int_{0}^{x}\left(e^{\frac{ye^{2 t}}{4}-\frac{t y}{2}+t} \left(e^t-1\right) \int_{0}^{t}e^{-\frac{ye^{2 u}}{4}+\frac{y u}{2}-u}d u \right)d t - x,\\
I^{1,1}(x,y)&=y\int_{0}^{x}\left(e^{-\frac{ty}{2}+ \frac{y e^{2 t}}{4}+t}\left(e^t-1\right)\int_{0}^{t}e^{-\frac{y e^{2 u}}{4}-u +\frac{y u}{2}}d u \right)d t,\\
I^{1,0}(x,y)&=y\int_0^x(e^t-1)(I^{0,0}(t,y)+I^{1,0}(t,y)-t-1)dt +\int_0^x (I^{1,0}(t,y)+I^{1,1}(t,y)+t)dt, 
\end{align*}
and
\begin{multline*}
I^{0,0}(x,y)=1+x+\frac{y}{2}\int_0^x\left((e^t-1)^2(I^{0,0}(t,y)+I^{1,0}(t,y)-t-1)\right)dt+\\
 +\int_0^x(I^{0,0}(t,y)+I^{1,0}(t,y)-t-1)dt + \int_0^x(I^{0,0}(t,y)+I^{0,1}(t,y)-t-1)dt.
 \end{multline*}
\end{theorem} 

\begin{proof} 
Let $\pi$ be a nonempty increasing permutation. We consider the decomposition of $\pi$ with respect to the  position of the entry 1.  From this 
\begin{align*} 
\I^{1,1}&=\left(\mathcal{V}^1\cup\{x\}\right) \star \{x\}^{\square}\star \left(\I^{0,1} \cup\I^{1,1}\cup\{x\}\right)\\
\I^{0,1}&=\left(\mathcal{V}^0\backslash\{\epsilon,x\}\right) \star \{x\}^{\square}\star \left(\I^{0,1} \cup\I^{1,1}\cup\{x\}\right)+\{x\}^{\square}\star \left(\I^{0,1} \cup\I^{1,1}\cup\{x\}\right)\\
\I^{1,0}&=\left(\mathcal{V}^1\cup\{x\}\right) \star \{x\}^{\square}\star \left(\I^{0,0} \cup\I^{1,0}\backslash\{\epsilon,x\}\right)+ \left(\I^{1,0} \cup\I^{1,1}\cup\{x\}\right)\star\{x\}^{\square}\\
\I^{0,0}&=\epsilon +\{x\}+ \left(\mathcal{V}^0\backslash\{\epsilon,x\}\right)\star \{x\}^{\square}\star\left(\I^{0,0}\cup\I^{1,0}\backslash\{\epsilon,x\}\right)+\{x\}^{\square}\star\left(\I^{0,0}\cup\I^{1,0}\backslash\{\epsilon,x\}\right)+\\&\qquad\qquad\qquad\qquad\qquad\qquad\qquad\qquad\qquad\qquad\left(\I^{0,0}\cup\I^{0,1}\backslash\{\epsilon,x\}\right)\star\{x\}^{\square},
\end{align*}
where $\mathcal{V}^a$, $a\in\{0,1\}$, is the set of valleyless permutations $\pi$ such that $\pi$ does not start with a descent if and only if $a=0$. From the symbolic method (cf. \cite{Flajolet2}) we obtain a system of functional equations which gives the desired result for $I^{a,b}(x,y)$, $a,b\in\{0,1\}$. We obtain the desired result for $I^\peak(x,y)$ by observing that: the number of peaks  and the number of valleys are equal in a permutation $\pi\in\I^{0,1}\cup\I^{1,0}$; the number of peaks equals to the number of valleys plus one in a permutation  $\pi\in\I^{0,0}\backslash\{\epsilon, 1\}$; and the number of peaks equals to the number of valleys minus one in a permutation  $\pi\in\I^{1,1}$. 
\end{proof}

We were unable to obtain solutions that could be presented in the paper. Then we leave it as an open question for the reader to attempt to find a simpler expression for the generating function $I^\peak(x,y)$.  

Let $\inc^\peak(n,k)$ the number of increasing permutations in $\Sn_n$ with exactly $k$ peaks. 
Then, we have the following \footnotesize
$$I^\peak:=[\inc^\peak(n,k)]_{n, k \geq 0}=\left(
\begin{array}{cccccc}
 1 & 0 & 0 & 0 & 0 & \cdots \\
 1 & 0 & 0 & 0 & 0 & \cdots \\
 2 & 0 & 0 & 0 & 0 & \cdots \\
 4 & 2 & 0 & 0 & 0 & \cdots \\
 8 & \framebox{\textbf{16}} & 0 & 0 & 0 & \cdots \\
 16 & 80 & 16 & 0 & 0 & \cdots \\
 32 & 341 & 211 & 0 & 0 & \cdots \\
 64 & 1361 & 1815 & 136 & 0 & \cdots \\
 128 & 5286 & 12988 & 3078 & 0 & \cdots \\
 \vdots & \vdots & \vdots& \vdots & \vdots&\ddots\\
\end{array}
\right).$$ \normalsize

For example, $\inc^\peak(4,1)=16$, the entry boxed in the array $I^\peak$ above, and the corresponding increasing permutations are
\begin{align*}
&1243, \ 1324, \ 1342,  \ 1423,  \ 1432, \ 2143, \ 2314, \ 2341,\\ &2413, \ 2431, \ 3142, \ 3241, \ 3412,  \ 3421, \ 4132,  \ 4231.
\end{align*}

\begin{coro}
The generating function for total number of peaks in all increasing  permutations of a given length is
\begin{align*}
\tot^\peak(x)&=\partial_y(I^\peak(x,y))\rvert_{y=1} \\
&=  2\frac{x^3}{3!} +16\frac{x^4}{4!}  + 112\frac{x^5}{5!}  + 763\frac{x^6}{6!}+ 5399\frac{x^7}{7!}+40496\frac{x^8}{8!} + O(x^{9}).
\end{align*}
\end{coro}
The sequence of coefficients of $\frac{x^n}{n!}$ in this series expansions does not appear in the OEIS.

\subsection{Alternating increasing permutations}

A permutation $\pi=\pi_1\pi_2\cdots \pi_n$ is \emph{alternating} if $\pi_1>\pi_2<\pi_3>\cdots$. The exponential  generating function for the number of alternating permutations is a classical result of Andr\'e  \cite{Andre}, and is given by 
$$\sum_{n\geq 0}E_n \frac{x^n}{n!}=\sec x + \tan x,$$
where $E_n$ corresponds to the Euler number (see \seqnum{A000111} in \cite{OEIS}).
A recent survey on alternating permutations was given by Stanley in \cite{Stanley2}.

\begin{theorem}
    The generating function $A(x)$ for the total number of increasing alternating permutations is 
\begin{align*}
    A(x)&=e^{\frac{x^2}{2}} (x+1)\int_{0}^xe^{-t^2}dt-e^{\frac{x^2}{2}}+2\\
&=1 + \frac{x}{1!}  + \frac{x^2}{2!} + 2\frac{x^3}{3!} + 5 \frac{x^4}{4!}  + 8 \frac{x^5}{5!}  +  33\frac{x^6}{6!} + 48\frac{x^7}{7!} +  279\frac{x^8}{8!}  + O(x^{8}).
\end{align*}
\end{theorem}
\begin{proof}
 Let $\mathcal{A}_{E}$ (resp. $\mathcal{A}_{O}$) be the set of increasing alternating permutations of even (resp. odd) length, respectively. From the decomposition with respect to the position of the entry 1, we obtain the symbolic equations:
 \begin{align*}
     \mathcal{A}_{O}&=\{x\} + 
\{x\} * \{x\}^{\square} *  \mathcal{A}_{O}, \\
\mathcal{A}_{E}&=\epsilon  +     \mathcal{A}_{O} *  \{x\}^{\square} + 
\{x\} * \{x\}^{\square} * (\mathcal{A}_{E}-\epsilon).
\end{align*}
Therefore, we obtain the differential  equations
 \begin{align*}
    \frac{d}{dx}A_{O}(x)&= 1 +  x A_{O}(x), \quad A_{O}(0)=0, \\
 \frac{d}{dx}A_{E}(x) &= A_{O}(x)  + x (A_E(x)-1), \quad A_E(0)=1.
\end{align*}
The solutions of these differential equations are
 \begin{align*}
  A_O(x)= e^{x^2/2} \int_{0}^xe^{-t^2/2}dt \quad \text{and} \quad A_E(x)= 2- e^{x^2/2} + x e^{x^2/2}  \int_{0}^xe^{-t^2/2}dt.
\end{align*}
Since $A(x)=A_O(x) + A_E(x)$, we obtain the desired result.
\end{proof}
Let $a_n$ denote the sequence of coefficients of $\frac{x^n}{n!}$ in $A(x)$.  Notice that $a_{2n+1}=(2n)!!$ and $a_{2n+2}=(2 n + 2)!! - (2 n + 1)!!$.

\section{Flattened Permutations}

In this section, we are interested in the bivariate  generating function
\begin{align*}
F^s(x,y)=\sum_{\pi \in\F}\frac{x^{|\pi|}}{|\pi|!}y^{\s(\pi)}=\sum_{n\geq 0}\frac{x^{n}}{n!}\sum_{\pi\in\F_n}y^{\s(\pi)},
\end{align*} 
where $\F$ is the set of flattened permutations, i.e. permutations consisting of runs arranged from left to right,  with the first entries of each run in increasing order. Here, the statistic  $\s$ corresponds to $\val$, $\run=1+\des=n-\asc=1+\peak$, and $\rlm$.
The coefficient of $\frac{x^ny^k}{n!}$ in $F^s(x,y)$ is the number of flattened permutations $\pi$ of  length $n$ satisfying $s(\pi)=k$. Note that the statistic $\run$ has already been investigate by Nabawanda et al.\ \cite{naraba} using combinatorial arguments and recursive formulas. 

\subsection{For the number of valleys: $\s=\val$} In this part, we focus on the enumeration of flattened permutations with respect to the length and the statistic $\s=\val$, which counts the number of valleys.

\begin{theorem}
The generating function for flattened permutations with respect to the
length and the number of valleys is
$$F^\val(x,y)=1+\int_{0}^{x}B^\val(t,y)dt,$$
where 
$$B^\val(x,y)=e^{ -xy+e^x y+x-y} \left(e^y \int _0^xe^{-e^{t} y+ty-t} \left(-e^{t}
   y+y+e^{t}-1\right)dt+1\right).$$
\end{theorem}
\begin{proof} We first define the subset $\mathcal{B}$ of permutations $\pi'\in \Sn$ having one of the following forms: ($i$) $\pi'=\epsilon$, or ($ii$) $\pi'=1(1+\pi_1)$ where $\pi_1\in \mathcal{B}$, or ($iii$) $\pi'=23\cdots k 1$ where $k\geq 2$, or ($iv$) $\pi=\pi_21\pi_1$ where $\texttt{red}(\pi_2)=12\cdots k$, $k\geq 1$ and $\texttt{red}(\pi_1)\in\mathcal{B}\backslash\{\epsilon\}$. Considering this definition, a nonempty flattened permutation $\pi$ is necessary of the form $\pi=1 (1+\pi')$ where $\pi'\in\mathcal{B}$. So, we first give the exponential generating function of the set $\mathcal{B}$ with respect to the length and the number of valleys.

From the above  decomposition of $\mathcal{B}$ we obtain the symbolic equation:
$$\B=\epsilon + \{x\}^{\square}\star \B + (\id-\epsilon) \star \{x\}^{\square} + (\id-\epsilon) \star \{x\}^{\square} \star(\B-\epsilon),$$
where the set $\id$ consists of all identity permutations of length $n\geq 0$. Therefore, we obtain the following functional equation:
$$B^\val(x,y)= 1 +  \int_{0}^x B^\val(t,y)dt + \int_{0}^x (e^t-1)dt + y\int_{0}^x (e^t -1)(B^\val(t,y)-1)dt.$$
Differentiating the previous equation, we obtain the following differential equation:
$$\frac{d}{dx} B^\val(x,y)=    B^\val(x,y) + e^x-1 + y (e^x -1) (B^\val(x,y)-1), \quad B^\val(0,y)=1.$$
The solution of this equation provides the expected expression of $B(x,y)$. Now, by considering the decomposition $\F=1+\{x\}^{\square}\star B$, we obtain the desired result. 
\end{proof}

We can use a symbolic software computation to obtain the first few terms of
the formal power series: 
\begin{multline*}
F^\val(x,y)=(2  x + e^x (x-1))+\frac{1}{4} \left(-2 e^x (x-2) x+e^{2 x} (2 x-3)+3\right) y \\+\frac{1}{36} \left(-9 e^{2 x} \left(2 x^2-2 x-7\right)+6 e^x \left(2 x^3-6 x^2-15\right)+2 e^{3 x}
   (6 x-11)+49\right) y^{2}+O(y^3).
   \end{multline*}

Let $f^\val(n,k)$ denote the number of flattened permutations in $\Sn_n$ with exactly $k$ valleys. The first few values of this sequence are
\footnotesize$$F^\val:=[f^\val(n,k)]_{n, k \geq 0}=\left(
\begin{array}{cccccc}
 1 & 0 & 0 & 0 & 0&\cdots \\
 1 & 0 & 0 & 0 & 0 &\cdots\\
 1 & 0 & 0 & 0 & 0&\cdots \\
 2 & 0 & 0 & 0 & 0&\cdots \\
 3 & 2 & 0 & 0 & 0&\cdots \\
 4 &\framebox{\textbf{11}} & 0 & 0 & 0&\cdots \\
 5 & 39 & 8 & 0 & 0 &\cdots\\
 6 & 114 & 83 & 0 & 0&\cdots \\
 7 & 300 & 522 & 48 & 0&\cdots \\
 8 & 741 & 2594 & 797 & 0&\cdots \\
 \vdots & \vdots      & \vdots       &  \vdots   & \vdots    &\ddots
\end{array}
\right).$$ \normalsize

For example, $f^\val(5,2)=11$, the entry boxed in the array $F^\val$ above, and the corresponding flattened permutations are
\begin{align*}
&12435, \quad 12534, \quad 13245, \quad 13254, \quad 13425, \quad
 13524,\\ & 14235, \quad 14253, \quad 14523,\quad
15234,  \quad 15243.
\end{align*}

\begin{coro}We have the following recursive formula for $n,k\geq 1$,
$f^\val(n+1,k)=b(n,k)$,
with 
$$b(n+1,k)= b(n,k)+\sum\limits_{i=1}^{n-1}{n \choose i}b(i,k-1)$$
with the initial conditions $b(0,0)=1$, $b(n,0)=n$ for $n\geq 1$, and $f^\val(0,0)=1$.
\end{coro}

\begin{coro}
The generating function for flattened  permutations  without valleys is
\begin{align*}
F^\val(x,0)&=2+x+(x-1)e^x\\
&=1  + \frac{x}{1!}  + \frac{x^2}{2!}  + 2\frac{x^3}{3!} +3 \frac{x^4}{4!}  + 4\frac{x^5}{5!}  + 5\frac{x^6}{6!} +6\frac{x^7}{7!}+ O(x^8),
\end{align*}
and the number of flattened permutations of length $n$ without valleys is $n-1$ for $n\geq 2$.
\end{coro}

\begin{coro}
The generating function  $F(x)$ for the total number of flattened  permutations is
\begin{align*}
   F(x)&:=F^{\val}(x,1)= 1 + \int_{0}^x e^{e^t-1}dt\\
   &=1  + \frac{x}{1!}  + \frac{x^2}{2!}  + 2\frac{x^3}{3!} +5\frac{x^4}{4!}  + 15\frac{x^5}{5!}  + 52\frac{x^6}{6!} +203\frac{x^7}{7!}+ O(x^8).
\end{align*}
\end{coro}
As expected (see \cite{manshawag}), the total number of flattened permutations of length $n\geq 0$ is given by the Bell number $B_{n-1}$ (sequence \seqnum{A000110}).

\begin{coro}
The generating function for the total number of valleys in all flattened permutations with respect to the length is given by  
\begin{align*}Tot^\val(x)&=\partial_y(F^\val(x,y))\rvert_{y=1}\\
&=\int _{0}^{x}\left(\!{\it Ei} \left(\left( 1,1 \right)-
{\it{Ei}}\left( 1,e^t \right)\right)e^{e^t} -(2+x+e^t)e^{e^t-1}\right){dt},
\end{align*}
where $$Ei(a,z)=\int_1^\infty e^{-tz}t^{-a}dt.$$
\end{coro}
The first few terms of the formal power series are
\begin{align*}Tot^\val(x)=2\frac{x^4}{4!}  + 11\frac{x^5}{5!}  + 55\frac{x^6}{6!} +280\frac{x^7}{7!}+1488\frac{x^8}{8!}+ O(x^9).
\end{align*}
The sequence of coefficients of $\frac{x^n}{n!}$ in this series expansions does not appear in \cite{OEIS}.

\subsection{For the number of runs: $\s=\run=1+\des=1+\peak=n-\asc$} 
In this part, we focus on the enumeration of flattened permutations with respect to the length and the statistic $\s=\run$ that counts the number of runs. A simple observation allows us to see that $\run=1+\des=1+\peak$ on the set of  flattened permutations.

The following theorem was first proved by Nabawanda et al.\ \cite{naraba}. We provide an alternative proof using the symbolic method.
\begin{theorem}
The generating function for flattened permutations with respect to the
length and the number of runs is
$$F^\run(x,y)= 1 + xy+ \int_{0}^x (B^\run(t,y)-1)dt,$$
where $$B^\run(x,y)=e^{-xy+x-y} \left(ye^{e^x y} +(1 -y)e^{x(y-1)+y}\right).
$$
\end{theorem}
\begin{proof} We make a similar proof as for Theorem 5.1 by considering the set $\mathcal{B}$ defined by the symbolic equation:
$$\B=\epsilon + \{x\}^{\square}\star \B + (\id-\epsilon) \star \{x\}^{\square} + (\id-\epsilon) \star \{x\}^{\square} \star(\B-\epsilon).$$
Thus we deduce 
$$B^\run(x,y)= 1 + xy+ \int_{0}^x (B^\run(t,y)-1)dt + y^2\int_{0}^x (e^t-1)dt + y\int_{0}^x (e^t -1)(B^\run(t,y)-1)dt.$$
Differentiating the previous equation, we obtain the following differential equation:
$$\frac{d}{dx} B^\run(x,y)=   y+ B^\run(x,y) -1+ y^2(e^x-1) + y(e^x -1) (B^\run(x,y)-1), \quad B^\run(0,y)=1.$$
The solution of this equation provides the expected expression of $B^\run(x,y)$. Now, by considering the decomposition $\F=1+\{x\}^{\square}\star B$, we obtain 
$$F^\run(x,y)= 1 + xy+ \int_{0}^x (B^\run(t,y)-1)dt,$$ which provides 
the desired result. 
\end{proof}

We can use a symbolic software computation to obtain the first few terms of
the formal power series 
\begin{multline*}
F^\run(x,y)=1+x + \left(e^x-1\right) y +\left(e^x (\sinh(x)-x)\right) y^{2}+\\\frac{1}{12} \left(6 e^x \left(x^2+1\right)-3 e^{2 x} (2 x+1)+2 e^{3 x}-5\right) y^{3} + O(y^4).
   \end{multline*}

Let $f^\run(n,k)$ denote the number of flattened permutations in $\Sn_n$ with exactly $k$ runs. The first few values of this sequence are 
\footnotesize $$F^\run:=[f^\run(n,k)]_{n, k \geq 0}=\left(
\begin{array}{cccccc}
 1 & 0 & 0 & 0 & 0&\cdots \\
 0 & 1 & 0 & 0 & 0 &\cdots\\
 0 & 1 & 0 & 0 & 0&\cdots \\
 0 & 1 & 1 & 0 & 0&\cdots \\
 0 & 1 & 4 & 0 & 0&\cdots \\
 0 &1&\framebox{\textbf{11}}  & 3 & 0&\cdots \\
 0 & 1 & 26 & 25 & 0 &\cdots\\
 0 & 1 & 57 & 130 & 15&\cdots \\
 0 & 1 & 120 & 546 & 210&\cdots \\
 \vdots & \vdots      & \vdots       &  \vdots   & \vdots    &\ddots
\end{array}
\right).$$ \normalsize

For example, $f^\run(5,2)=11$, the entry boxed in the array $F^\run$ above, and the corresponding flattened permutations are
\begin{align*}
&12354, \quad 12435, \quad 12453, \quad 12534, \quad 13245, \quad 13425, \\\quad
 &13452, \quad 13524, \quad 14235, \quad 14523,\quad
15234.
\end{align*}

Nabawanda et al.\ \cite{naraba} proved the following recurrence relation
$$f^\run(n,k)=kf^\run(n-1,k) + (n-2)f^\run(n-2,k-1).$$

\begin{coro}
The generating function for the total number of runs in all flattened permutations with respect to the length is given by  
\begin{align*}Tot^\run(x)&=\partial_y(F^\val(x,y))\rvert_{y=1}=x +\int_{0}^{x}\left(-e^{e^t-1} t+e^{t+e^t-1}-1\right)d t.\\
&=x  +\frac{x^2}{2!}  +3\frac{x^3}{3!}  +9\frac{x^4}{4!}  + 32\frac{x^5}{5!}  + 128\frac{x^6}{6!} +565\frac{x^7}{7!}+2719\frac{x^8}{8!} + O(x^9).
\end{align*}
\end{coro}

Let $\F^{\val=\des}_n$ be the set of flattened permutations $\pi$ of length $n$ such that   $\val(\pi)=\des(\pi)$. We set $\F^{\val=\des}=\sum_{n\geq 0}\F^{\val=\des}_n$. 

\begin{theorem}  The cardinality of $\F^{\val=\des}_n$ is given by the $n$-th Gould number (see \seqnum{A040027} in \cite{OEIS}):
$$|\F^{\val=\des}_n|=\sum\limits_{k=0}^{n-2}\sum\limits_{j=k}^{n-2}\stirling{j}{k}k^{n-2-j},$$
where $\stirling{j}{k}$ is the Stirling number of second kind.
\end{theorem}
\begin{proof} The proof is obtained {\it mutatis mutandis} as for Theorem 2.8.
\end{proof}

\subsection{For the number of right-to-left minima: $\s=\rlm$}

In this part, we focus on the enumeration of flattened permutations with respect to the length and the statistic $\s=\rlm$ that counts the number of right-to-left minima.

\begin{theorem}
The generating function for flattened permutations with respect to the
length and the number of right-to-left minima is
\begin{align*}F^\rlm(x,y)&= 1 + y\int_{0}^x B^\rlm(t,y)dt,
\end{align*}
where $B^\rlm(x,y)=e^{y \left(e^{x}-1\right)}$. \end{theorem}
\begin{proof} We make a similar proof as for Theorem 5.1 and Theorem 5.6 by considering the set $\mathcal{B}$ defined by the symbolic equation:
$$\B=\epsilon + \{x\}^{\square}\star \B + (\id-\epsilon) \star \{x\}^{\square} + (\id-\epsilon) \star \{x\}^{\square} \star(\B-\epsilon).$$
Thus we deduce the functional equation 
$$B^\rlm(x,y)= 1 + y\int_{0}^x B^\rlm(t,y)dt + y\int_{0}^x (e^t-1)dt + y\int_{0}^t (e^t -1)(B^\rlm(t,y)-1)dt.$$
Differentiating the previous equation, we obtain the following differential equation:
$$\frac{d}{dx} B^\rlm(x,y)=  yB^\rlm(x,y)+ y(e^x-1) + y(e^x -1) (B^\rlm(x,y)-1), \quad B^\rlm(0,y)=1.$$
The solution of this equation provides the expected expression of $B^\rlm(x,y)$. Now, by considering the decomposition $\F=1+\{x\}^{\square}\star B$, we obtain 
$$F^\rlm(x,y)= 1 + y\int_{0}^x B^\rlm(t,y)dt,$$ which provides 
the desired result. 
\end{proof}

The first few terms of the formal power series are
\begin{multline*}
F^\rlm(x,y)=1+yx+ \left( e^x-x-1 \right) {y}^{2}+ \frac{1}{4} \left(2 x-4 e^x+e^{2 x}+3\right) {y}^{3}+ \\\frac{1}{36} \left(-6 x+18 e^x-9 e^{2 x}+2 e^{3 x}-11\right) {y}^{4} + O(y^5).
   \end{multline*}

Let $f^\rlm(n,k)$ the number of flattened permutations in $\Sn_n$ with exactly $k$ right-to-left minima. The first few values of this sequence are
\footnotesize
$$F^\rlm:=[f^\rlm(n,k)]_{n, k \geq 0}=\left(
\begin{array}{cccccccc}
 1 & 0 & 0 & 0 & 0&0&0&\cdots \\
 0 & 1 & 0 & 0 & 0 &0&0&\cdots\\
 0 & 0 & 1 & 0 & 0&0&0&\cdots \\
 0 & 0 & 1 & 1 & 0&0&0&\cdots \\
 0 & 0 & 1 & 3 & 1&0&0&\cdots \\
 0 &0&  1& \framebox{\textbf{7}} & 6&1&0&\cdots \\
 0 & 0 & 1 & 15 & 25 &10&1&\cdots\\
 0 & 0 & 1 & 31 & 90&65&15&\cdots \\
 \vdots & \vdots      & \vdots       &  \vdots   & \vdots & \vdots  & \vdots  &\ddots
\end{array}
\right).$$\normalsize
For example, $f^\rlm(5,3)=7$, the entry boxed in the array $F^\rlm$ above, and the corresponding flattened permutations are
\begin{align*}
 12453, \quad 13254, \quad 13425, \quad 13524, \quad 14253,\quad 14523,\quad 15243. 
\end{align*}

\begin{coro}  For $n,k\geq 1$, we have
$$f^\rlm(n,k)=\stirling{n-1}{k-1}.$$
\end{coro}

Notice that $f^\rlm(n,k)$ equals the Stirling number of the second kind $\stirling{n-1}{k-1}$ that also counts the number of set partitions of $[n-1]$ with $k-1$ blocks. We decompose a flattened permutation  $\pi=a_1\pi_2a_2 \ldots \pi_ka_k$ where $a_1=1, a_2, \cdots, a_k$ are the right-to-left minima of $\pi$, and we construct the partition $B=B_1/B2/\cdots /B_k$ of $[n-1]$,  where $B_i=\texttt{dom}((\pi_i-1)(a_i-1))$ for $2\leq i\leq k$.

\begin{coro}
The generating function for the total number of right-to-left minima in all flattened permutations with respect to the length is given by $n$-th Bell number.
\end{coro}

\end{document}